\documentclass[11pt, leqno]{article}

\usepackage[T1]{fontenc}
\usepackage[utf8]{inputenc}
\usepackage[english]{babel}

\usepackage{amsmath}
\usepackage{geometry}
\usepackage{color}
\usepackage{xcolor}
\usepackage{enumitem}
\usepackage[linktoc=all]{hyperref}

\usepackage{graphicx}
\usepackage{wrapfig}
\usepackage{caption}
\usepackage{float}
\usepackage{tikz}

\usepackage{mathtools}
\usepackage{amssymb}
\usepackage{wasysym}
\usepackage{textcomp}
\usepackage{stmaryrd}
\usepackage{esint}

\usepackage{amsthm}
\usepackage{thmtools}

\usepackage{xspace}

\allowdisplaybreaks[1]

\colorlet{myred}{red!50!black}
\colorlet{mylightblue}{blue!50!black}
\colorlet{mydarkblue}{blue!80!black}
\colorlet{mygreen}{green!50!black}

\hypersetup{%
	colorlinks = true,
	linkcolor = myred,
	citecolor = mylightblue,
	urlcolor = mydarkblue
}

\usetikzlibrary{matrix}
\usetikzlibrary{arrows,decorations.markings}
\usetikzlibrary{intersections}
\tikzstyle{fleche}=[->, >=latex]
\tikzset{
	xa/.store in=\xa, xa/.default=0, xa=0,%
	xb/.store in=\xb, xb/.default=0, xb=0,%
	xc/.store in=\xc, xc/.default=0, xc=0,%
}

\newcommand*{\N}{\mathbb{N}}

\newcommand*{\R}{\mathbb{R}}
\newcommand*{\C}{\mathbb{C}}
\newcommand*{\dd}{\mathrm{d}}

\newcommand*{\deq}{\overset{\scriptscriptstyle{\mathrm{def}}}{=}}

\newcommand{\ri}{\mathrm{i}}
\newcommand{\re}{\mathrm{e}}
\newcommand{\de}{\mathrm{d}}

\newcommand*{\Norm}[1]{\left\lVert#1\right\rVert}
\newcommand*{\norm}[1]{\lVert#1\rVert}

\newcommand*{\absolute}[1]{\lvert#1\rvert}

\newcommand*{\scalp}[1]{\langle #1\rangle}

\newcommand*{\Cal}[1]{\mathcal{#1}}

\newcommand{\rphi}{\phi_\text{re}}
\newcommand{\iphi}{\phi_\text{im}}

\newcommand{\diffusion}{d}

\newtheorem{theorem}{Theorem}[section]

\newtheorem{proposition}[theorem]{Proposition}
\newtheorem{lemma}[theorem]{Lemma}

\theoremstyle{definition}

\theoremstyle{remark}
\newtheorem{remark}[theorem]{Remark}

\title{Global existence and decay of small solutions in a viscous half Klein-Gordon equation}

\author{Louis Gar\'enaux$^*$ and Bj\"orn de Rijk\thanks{Institute for Analysis, Karlsruhe Institute of Technology, Englerstra\ss e 2, 76131 Karlsruhe, Germany}}

\date{December 27, 2022}

\begin{document}

\maketitle


\begin{abstract}
We establish global existence and decay of solutions of a viscous half Klein-Gordon equation with a quadratic nonlinearity considering initial data, whose Fourier transform is small in $L^1(\R) \cap L^\infty(\R)$. Our analysis relies on the observation that nonresonant dispersive effects yield a transformation of the quadratic nonlinearity into a subcritical nonlocal quartic one, which can be controlled by the linear diffusive dynamics through a standard $L^1$-$L^\infty$-argument. This transformation can be realized by applying the normal form method of Shatah or, equivalently, through integration by parts in time in the associated Duhamel formula.

\end{abstract}

\section{Introduction}

We establish global existence and decay of solutions with small initial data in the viscous half Klein-Gordon equation
\begin{equation}
\label{e:intro:main-eq}
\partial_t u - \lambda(\partial_x) u = B(u, u),
\end{equation}
with $t \geq 0$, $x \in \R$ and $u(t,x) \in \C$, where the pseudo-differential operator $\lambda$ is defined through its Fourier symbol
\begin{equation} \label{e:Fourier}
\hat{\lambda}(k) \deq - d k^2 - \ri\scalp{k},
\end{equation}
with $d > 0$ representing viscosity, and where $B$ is a bilinear form given by
\begin{equation} \label{e:defB}
B(u_1, u_2) = \Cal{F}^{-1}\left(k \mapsto \int_\R \hat{B}(k, l) \hat{u}_1(k-l) \hat{u}_2(l) \dd l\right),
\end{equation}
with $\hat{B} \in L^\infty(\R^2)$. Here, we use the Japanese bracket notation $\scalp{k} \deq \sqrt{1+k^2}$ and $\Cal{F}$ stands for the Fourier transform
\begin{align*} \hat{u}(k) = (\Cal{F}u)(k) = \int_{\R} \re^{-\ri kx} u(x) \de x.\end{align*} 
We remark that the case $\hat{B}(k,l) = \smash{\frac{1}{2\pi}}$ corresponds to the standard quadratic nonlinearity $B(u, u) = u^2$.

Our study of~\eqref{e:intro:main-eq} is motivated by the one-dimensional Klein-Gordon equation with viscoelastic dissipation~\cite{AVI,AVR,ABIK,SHIB}, given by
\begin{align}
u_{tt} + \beta^2\left(\gamma - \partial_x^2\right) u - \alpha \partial_x^2 u_t = N(u), \label{e:vKG}
\end{align}
with $x \in \R$, $t \geq 0$, $u(x,t) \in \R$ and parameters $\alpha,\beta \geq 0$ and $\gamma > 0$, where $N \in C^1(\R,\R)$ is a nonlinearity with $N(0), N'(0) = 0$. After rewriting~\eqref{e:vKG} as a nonlinear evolution system, the Fourier symbol of its linearization about the trivial state possesses the eigenvalues 
\begin{align}\label{e:eigenvalues}
 \lambda_{\pm}(k) = -\dfrac{1}{2}\alpha k^2 \pm \sqrt{\dfrac{1}{4}\alpha^2k^4 - \beta^2(\gamma+k^2)},
\end{align}
which correspond to the Fourier symbol of the pseudo-differential operator $\lambda$ in~\eqref{e:intro:main-eq} up to the term $\frac{1}{4}\alpha^2k^4$ upon setting $\alpha = 2d$, $\beta = 1$ and $\gamma = 1$. We note, as long as $\beta > 0$, it is always possible to realize $\beta = 1$ and $\gamma = 1$ by rescaling time, space and the parameter $\alpha$ in~\eqref{e:vKG}.

It is well-known that the asymptotics of solutions with small initial data in (damped) nonlinear wave equations such as~\eqref{e:vKG} is delicate. More precisely, the temporal decay exhibited by the linear terms in~\eqref{e:vKG} is algebraic and at rate $\smash{t^{-\frac{1}{2}}}$ in $L^\infty(\R)$ for sufficiently localized initial data, cf.~\cite{ABIK,SHIB}. For nonlinearities $N(u) = u^q$ and functions $u(t)$ decaying at such rate in $L^\infty(\R)$ we thus obtain $\smash{\|N(u(t))\|_{L^p} \leq t^{-\frac{q-1}{2}} \|u(t)\|_{L^p}}$ for $1 \leq p \leq \infty$ and $q \in \N$, showing that in case $q \leq 3$ we are at or below the threshold of integrability. Hence, standard arguments providing control on quadratic or cubic nonlinear terms in~\eqref{e:vKG} by the linear dynamics through iterative estimates on the Duhamel formulation do not apply and the asymptotics of solutions with small initial data could depend on the precise structure of the nonlinearity.

In fact, in the purely dispersive case with $\alpha = 0$, i.e.~without viscoelastic dissipation, all solutions of~\eqref{e:vKG} with small, compactly supported initial data blow up in finite time for $N(u) = |u|u$, cf.~\cite{KETA}. Moreover, also in the purely dissipative case obtained by setting $\beta = 0$, one readily observes that solutions of~\eqref{e:vKG} with (arbitrarily small) Gaussian initial data $u(x,0) = \eta \smash{\re^{-x^2}}$, $u_t(x,0) = 0$ with $\eta > 0$ blow up in finite time for $N(u) = u^2$.

On the other hand, there are various global existence results for solutions of~\eqref{e:vKG} with small initial data in the purely dispersive case with $\alpha = 0$ and $\beta > 0$, allowing for nonlinear terms of the form $N(u) = a_0 u^2 + a_1 u^3$ with $a_0, a_1 \in \R$, see~\cite{DELO,HANA3,HANA,STIN} and references therein. These results rely on the classical observation of Shatah~\cite{SHAT} that these nonlinear terms satisfy a nonresonance condition in case of nonzero dispersion and can be transformed into nonlocal higher-order terms using normal forms.

Although global existence of solutions of~\eqref{e:vKG} in the case with viscoelastic dissipation, i.e.~in case $\alpha, \beta > 0$, has been established in~\cite{AVI} for higer spatial dimensions, where the nonlinearities are subcritical, such results seem to be unavailable in spatial dimension one. In this contribution we make a first step in this direction by proving global existence and decay of solutions of~\eqref{e:intro:main-eq} with initial data whose Fourier transform is small in $L^1(\R) \cap L^\infty(\R)$. That is, we prove the following result.

\begin{theorem}
\label{t:main-result}
There exist $\delta>0$ and $C>0$ such that, if the Fourier transform of $u_0 \in L^2(\R)$ satisfies $\hat{u}_0 \in L^1(\R) \cap L^\infty(\R)$ with
\begin{equation*}
E_0 \deq \norm{\hat{u}_0}_{L^1\cap L^\infty} \leq \delta,
\end{equation*}
then there exists a unique global mild solution $u\in C\big([0,\infty),L^2(\R)\big)$ of~\eqref{e:intro:main-eq} with initial condition $u_0$. Furthermore, the solution $u(t)$ satisfies $\hat{u}(t) \in L^1(\R) \cap L^\infty(\R)$ and obeys the estimates
\begin{equation*}
\norm{u(t)}_{L^\infty} \leq \norm{\hat{u}(t)}_{L^1} \leq C E_0(1+t)^{-\frac{1}{2}}, 
\hspace{4em}
\norm{u(t)}_{L^2} \leq C E_0(1+t)^{-\frac{1}{4}}, 
\end{equation*}
\begin{equation*}
\norm{\hat{u}(t)}_{L^\infty} \leq C E_0, 
\end{equation*}
for all $t \geq 0$.
\end{theorem}

As explained before, the quadratic nonlinearity in~\eqref{e:intro:main-eq} cannot be controlled by the linear dynamics through standard iterative estimates on the Duhamel formulation. Thus, inspired by the analyses in~\cite{DELO,HANA3,HANA,STIN}, we show that the quadratic nonlinearity in~\eqref{e:intro:main-eq} satisfies a nonresonance condition and can be transformed into a subcritical quartic one using the normal form method of Shatah. We note that the absence of resonances strongly hinges on the dispersive character of~\eqref{e:intro:main-eq}. That is, without the dispersive term $\ri \scalp{k}$ in the Fourier symbol of the pseudo-differential operator $\lambda$ in~\eqref{e:intro:main-eq} Theorem~\ref{t:main-result} fails, see Remark~\ref{rem:disp}. 

Our analysis confirms that the normal form method works well in the presence of dissipative terms. In fact, after applying the normal form transform, the remaining quartic nonlinearity can be controlled by the linear diffusive term $d\partial_x^2$ in~\eqref{e:intro:main-eq} through standard iterative $L^1$-$L^\infty$-estimates, cf.~\cite[Chapter~14.1.3]{SUBOOK}, on the associated Duhamel formulation. This allows us to handle a significantly larger (in the sense of regularity and localization) class of initial data compared to the aforementioned works~\cite{DELO,HANA3,HANA,STIN} treating the purely dispersive case, where initial conditions $u_0$ are at least two times weakly differentiable and more strongly localized in the sense that at least $x \mapsto \scalp{x} u_0(x)$ must lie in $L^2(\R)$.

The normal form method is equivalent to an integration-by-parts procedure with respect to time in the Duhamel formulation. Thereby, it relates to the space-time resonances method~\cite{GERM,GMS,GMS2} as developed by Germain, Masmoudi and Shatah. This method is designed to capture dispersive decay in Duhamel-based arguments by blending (non$\text{-}$)stationary phase theory with bi- and multilinear analysis of oscillatory integrals arising in Fourier space. We refer to~\cite{GERMAIN,LAN,SHAT2} for more details and background. We explore the connection to the space-time resonances method in our setting by showing that Theorem~\ref{t:main-result} can also be proved through integration by parts with respect to time in the Duhamel formula associated with~\eqref{e:intro:main-eq}. 

Our paper is structured as follows. First, we employ the normal form method  in~\S\ref{sec:normalform} to prove Theorem~\ref{t:main-result}. Next, we relate in~\S\ref{sec:spacetime} to the space-time resonances method by showing how Theorem~\ref{t:main-result} follows through integration by parts with respect to time in the Duhamel formulation. Finally, Appendix~\ref{app} contains the (technical) analysis of the phase function, whose zeros correspond to possible resonances of the quadratic nonlinearity in~\eqref{e:intro:main-eq}.

\begin{remark}\label{rem:disp}
The necessity of dispersion becomes apparent by setting $\lambda(\partial_x) = d\partial_x^2$ and $B(u,u) = u^2$ in~\eqref{e:intro:main-eq}, which yields the heat equation with quadratic nonlinearity in which all nonnegative nontrivial initial data blow up in finite time~\cite{FUJI}.
\end{remark}

\begin{remark}\label{rem:KG}
For small frequencies we expect that our analysis can be transferred from the viscous half Klein-Gordon equation~\eqref{e:intro:main-eq} to the ``full'' Klein-Gordon equation~\eqref{e:vKG} with viscoelastic dissipation. Yet, different arguments are necessary for large frequencies, since the dispersive character of the Fourier symbol associated with the first-order formulation of~\eqref{e:vKG} disappears with its eigenvalues~\eqref{e:eigenvalues} becoming real for large $k$. In contrast, the Fourier symbol~\eqref{e:Fourier} of the pseudo-differential operator $\lambda$ in~\eqref{e:intro:main-eq} maintains its nonzero imaginary part for large $k$. 
\end{remark}

\textbf{Notation.}  Let $S$ be a set, and let $A,B\colon S\to \R$. Throughout the paper, the expression ``$A(x) \lesssim B(x)$ for $x \in S$'', means that there exists a constant $C > 0$, independent of $x$, such that $A(x) \leq CB(x)$ holds for all $x \in S$.

\textbf{Acknowledgement.} Funded by the Deutsche Forschungsgemeinschaft (DFG, German Research Foundation) -- Project-ID 258734477 -- SFB 1173.
 

\section{Normal form transform} \label{sec:normalform}

In this section we prove Theorem~\ref{t:main-result} using the normal form method of Shatah, cf.~\cite{DELO,HANA3,HANA,SHAT,STIN}. 

\subsection{Change of variable}

In the following proposition we introduce the relevant change of variables that transforms the quadratic nonlinearity in~\eqref{e:intro:main-eq} into a quartic one. We will work in the Banach space
\begin{align*} X \deq \left\{u \in L^2(\R) : \hat{u} \in L^1(\R) \cap L^\infty(\R)\right\},\end{align*}
equipped with the norm $\|u\|_X = \|\hat{u}\|_{L^1 \cap L^\infty}$. The space $X$ can be recognized as the inverse image of $L^1(\R) \cap L^\infty(\R)$ under the Fourier transform ${\cal F} \colon L^2(\R) \to L^2(\R)$.

\begin{proposition}
\label{p:normal-form}
Let $1 \leq p \leq \infty$. There exist a constant $C > 0$ and bounded bilinear and a trilinear forms $A_2$ and $A_3$ on $X$ obeying the estimates
\begin{align}
\norm{\Cal{F}(A_2(u,u))}_{L^p} & {} \leq C \norm{\hat{u}}_{L^1} \norm{\hat{u}}_{L^p}, \label{e:A2est}\\ 
\norm{\Cal{F}(A_3(u,u,u))}_{L^p} & {} \leq C \norm{\hat{u}}_{L^1}^2 \norm{\hat{u}}_{L^p}, \label{e:A3est}
\end{align}
such that, if $u \in C\big([0,\infty),X\big) \cap C^1\big((0,\infty),X\big)$ is a classical solution of~\eqref{e:intro:main-eq}, then $w \in  C\big([0,\infty),X\big) \cap C^1\big((0,\infty),X\big)$ given by
\begin{equation}
\label{e:normal-form:def-w}
w = u + A_2(u,u) + A_3(u,u,u)
\end{equation}
is a classical solution of
\begin{equation}
\label{e:weq}
\partial_t w - \lambda(\partial_x) w = Q(u),
\end{equation}
where $Q \colon X \to X$ is the quartic nonlinearity given by 
\begin{equation}
\label{e:def-Q}
\Cal{F}(Q(u))(k) = \int_{\R^3} \hat{Q}(k,l,m,n) \hat{u}(k-l) \hat{u}(l-m) \hat{u}(m-n) \hat{u}(n) \dd n \dd m \dd l.
\end{equation}
with $\hat{Q} \in L^\infty(\R^4)$, satisfying
\begin{align} \label{e:Qest}
\norm{\Cal{F}(Q(u))}_{L^p} & {} \leq C \norm{\hat{u}}_{L^1}^3 \norm{\hat{u}}_{L^p}.
\end{align}
\end{proposition}
\begin{proof}
We begin with the simpler ansatz $w = u + A_2(u,u)$, where $A_2$ is a bilinear form on $X$ defined by its action on Fourier space as
\begin{equation*}
\Cal{F}(A_2(u_1, u_2))(k) = \int_\R \hat{A}_2(k, l) \hat{u}_1(k-l) \hat{u}_2(l) \dd l,
\end{equation*}
and $\hat{A}_2 \in L^\infty(\R^2)$ is a function we are about to determine. Using the fact that $A_2$ is bilinear, we compute $\partial_t A_2(u,u) = A_2(\partial_t u, u) + A_2(u, \partial_t u)$ for $t > 0$. Thus, using equation~\eqref{e:intro:main-eq}, we arrive at
\begin{equation*}
\partial_t w - \lambda(\partial_x) w = B(u,u) - [\lambda(\partial_x), A_2](u,u) + T(u),
\end{equation*}
for $t > 0$, where $T(u)$ is the trilinear term given by $T(u) \deq A_2(B(u,u),u) + A_2(u, B(u,u))$, and where the square brackets stand for the linear-bilinear commutator
\begin{equation*}
[\lambda(\partial_x), A_2](u_1,u_2) \deq \lambda(\partial_x) A_2(u_1, u_2) - A_2(\lambda(\partial_x)u_1, u_2) - A_2(u_1, \lambda(\partial_x)u_2).
\end{equation*}
We see here that imposing the condition $[\lambda(\partial_x), A_2] = B$ allows us to replace quadratic by cubic terms in the equation satisfied by $w$. One readily verifies
\begin{equation*}
\Cal{F}\big([\lambda(\partial_x), A_2](u_1, u_2)\big)(k) = -\int_\R \phi(k,l)\hat{A}_2(k,l) \hat{u}_1(k-l) \hat{u}_2(l) \dd l,
\end{equation*}
where $\phi \colon \R^2 \to \C$ is given by 
\begin{equation*}
\phi(k,l) \deq -\hat{\lambda}(k) + \hat{\lambda}(k-l) + \hat{\lambda}(l).
\end{equation*}
Thus, setting $\hat{A}_2 \in L^\infty(\R^2)$ equal to $\hat{A}_2(k,l) = \frac{\hat{B}(k,l)}{\phi(k,l)}$, which is well-defined by Proposition~\ref{p:phi-study:main-result}, ensures that $[\lambda(\partial_x), A_2] = B$. This leads to 
\begin{equation*}
\partial_t w - \lambda(\partial_x) w = T(u),
\end{equation*}
for $t > 0$. Moreover, we estimate
\begin{equation*}
\norm{\Cal{F}(A_2(u,u))}_{L^p} \leq \Norm{\int_\R \norm{\hat{A}_2}_{L^\infty(\R^2)} \absolute{\hat{u}(k-l)} \absolute{\hat{u}(l)} \dd l}_{L^p} = \norm{\hat{A}_2}_{L^\infty(\R^2)} \Norm{\big.\absolute{\hat{u}} * \absolute{\hat{u}}}_{L^p}.
\end{equation*}
Applying Young's convolution inequality then leads to the estimate~\eqref{e:A2est}. 

Now, we work with the more advanced change of variable~\eqref{e:normal-form:def-w}, where the trilinear form $A_3$ on $X$ is defined by its action on Fourier space as
\begin{equation*}
\Cal{F}(A_3(u_1, u_2, u_3))(k) = \int_{\R^2} \hat{A}_3(k, l, m) \hat{u}_1(k-l) \hat{u}_2(l-m) \hat{u}_3(m) \dd m \dd l,
\end{equation*}
where $\hat{A}_3 \in L^\infty(\R^3)$ is to be determined. We now essentially follow the same steps as above, keeping the previous definition of $A_2$. Inserting the ansatz~\eqref{e:normal-form:def-w} into~\eqref{e:intro:main-eq}, we obtain 
\begin{equation*}
\partial_t w - \lambda(\partial_x) w = T(u) - [\lambda(\partial_x), A_3](u,u,u) + Q(u),
\end{equation*}
for $t > 0$, where the square brackets stand for the linear-trilinear commutator
\begin{align*}
[\lambda(\partial_x), A](u_1,u_2,u_3) &\deq \lambda(\partial_x) A(u_1,u_2,u_3) - A(\lambda(\partial_x) u_1,u_2,u_3) - A(u_1,\lambda(\partial_x) u_2,u_3)\\
&\qquad - \, A(u_1,u_2,\lambda(\partial_x) u_3),
\end{align*}
and where the quartic term $Q(u)$ is given by 
$$Q(u) \deq A_3(B(u,u),u, u) + A_3(u, B(u,u),u) + A_3(u, u, B(u,u)),$$ 
which corresponds to expression~\eqref{e:def-Q} with 
\begin{align}
\label{e:def-Q-hat}
\begin{split}
\hat{Q}(k,l,m,n) \deq {} & \hat{A}_3(k,m,n)\hat{B}(k-m,l-m) + \hat{A}_3(k,l,n)\hat{B}(l-n,m-n)\\
& {} + \hat{A}_3(k,l,m)\hat{B}(m,n).
\end{split}
\end{align}
One readily computes
\begin{equation*}
\Cal{F}\big([\lambda(\partial_x), A_3](u_1, u_2, u_3)\big)(k) = \int_{\R^2} \psi(k,l,m) \hat{A}_3(k,l,m) \hat{u}_1(k-l) \hat{u}_2(l-m) \hat{u}_3(m) \dd m \dd l,
\end{equation*}
where $\psi \colon \R^3 \to \C$ is given by
\begin{equation*}
\psi(k,l,m) \deq \hat{\lambda}(k) - \hat{\lambda}(k-l) - \hat{\lambda}(l-m) - \lambda(m) = -\phi(k,l) - \phi(l,m).
\end{equation*}
On the other hand, we have
\begin{align*}
&\Cal{F}(T(u))(k) = \int_{\R^2} \left(\hat{A}_2(k,l) \hat{B}(k-l, m) \hat{u}(k-l-m)\hat{u}(m)\hat{u}(l)\right.\\
&\qquad\qquad\qquad\qquad \left.\phantom{\int_{\R^2}} + \hat{A}_2(k,l) \hat{B}(l, m) \hat{u}(k-l)\hat{u}(l-m)\hat{u}(m) \right) \dd m \dd l, \\
&= \int_{\R^2} \left(\hat{A}_2(k,m) \hat{B}(k-m, l-m) +  \hat{A}_2(k,l) \hat{B}(l, m) \right) \hat{u}(k-l)\hat{u}(l-m)\hat{u}(m) \dd m \dd l.
\end{align*}
Thus, setting $\hat{A}_3 \in L^\infty(\R^3)$ equal to
\begin{equation}
\label{e:def-A3-hat}
\hat{A}_3(k,l,m) = \frac{1}{\phi(k,l) + \phi(l,m)} \left(\frac{\hat{B}(k,m) \hat{B}(k-m, l-m)}{\phi(k,m)} + \frac{\hat{B}(k,l) \hat{B}(l,m)}{\phi(k,l)}\right)
\end{equation}
which is well-defined by Proposition~\ref{p:phi-study:main-result}, ensures that $[\lambda(\partial_x), A_3](u,u,u) = T(u)$. Because $\hat{A}_3$ and $\hat{B}$ are in $L^\infty(\R^3)$ and $L^\infty(\R^2)$, respectively, we find that $\hat{Q}$ lies in $L^\infty(\R^4)$. 

Finally, the proof of the inequalities~\eqref{e:A3est} and~\eqref{e:Qest} goes along the same line as the estimate~\eqref{e:A2est}: we simply observe that $\hat{A}_3$ and $\hat{Q}$ are in $L^\infty(\R^3)$ and $L^\infty(\R^4)$, respectively, recall their expressions~\eqref{e:def-Q-hat} and~\eqref{e:def-A3-hat}, and a direct iteration of Young's inequality for convolution products leads to a constant $C_n > 0$ such that 
\begin{equation*}
\norm{u_1 * \dots * u_n}_{L^p} \leq C_n \norm{u_1}_{L^{p_1}} \cdots \norm{u_n}_{L^{p_n}}, 
\hspace{4em}
(n-1) + \frac{1}{p} = \sum_{i=1}^n \frac{1}{p_i}, 
\end{equation*}
which we apply with $p_1 = p$ and $p_i=1$ for $i = 2,\ldots,n$.
\end{proof}

\subsection{Linear bound}

We establish the relevant bound on the linear semigroup generated by the Fourier multiplier $\hat{\lambda}$ defined in~\eqref{e:Fourier}.

\begin{lemma}
\label{l:linear-bounds}
Let $1 \leq p \leq \infty$. There exists $C > 0$ such that we have 
\begin{equation}
\label{e:linear-bound}
\norm{k \mapsto \re^{t\hat{\lambda}(k)} \hat{u}_0(k)}_{L^p} \leq C \frac{\norm{\hat{u}_0}_{L^{p} \cap L^\infty}}{(1+t)^{\frac{1}{2p}}}.
\end{equation}
for all $t \geq 0$.
\end{lemma}
\begin{proof}
Let $1 \leq q, \alpha \leq \infty$ be such that $\frac{1}{p} = \frac{1}{q} + \frac{1}{\alpha}$. There exists a constant $C > 0$ such that
\begin{equation*}
\norm{k\mapsto \re^{t\hat{\lambda}(k)} \hat{u}_0(k)}_{L^p} = \norm{k\mapsto \re^{-tdk^2} \hat{u}_0(k)}_{L^p} \leq \norm{\hat{u}_0}_{L^q} \norm{k\mapsto \re^{-tdk^2}}_{L^\alpha} \leq C\norm{\hat{u}_0}_{L^q} t^{-\frac{1}{2\alpha}},
\end{equation*}
for $t > 0$, where we used H\"older's inequality and the standard integral identity
\begin{align} \int_\R \re^{-\diffusion k^2 t} \de k = \frac{\sqrt{\pi}}{\sqrt{\diffusion t}}, \qquad t > 0.\label{e:intid}\end{align} 
To obtain~\eqref{e:linear-bound} we combine the above estimate with $p = q$ and $\alpha = \infty$ for small times $t\in (0,1)$ with the estimate with $q = \infty$ and $\alpha = p$ for large times $t\geq 1$. 
\end{proof}

\subsection{Closing the nonlinear argument}

We prove Theorem~\ref{t:main-result} by first applying the change of variables~\eqref{e:normal-form:def-w} which transforms~\eqref{e:intro:main-eq} into its ``normal form''~\eqref{e:weq}. Subsequently, we close a nonlinear argument by iteratively estimating the right-hand side of the Duhamel formulation associated with~\eqref{e:weq} in $X$ using the linear and multilinear bounds, established in Lemma~\ref{l:linear-bounds} and Proposition~\ref{p:normal-form}, respectively.

\begin{proof}[Proof -- Theorem~\ref{t:main-result}]
Let $1 \leq p \leq \infty$ and $u_0 \in X$. By standard local existence theory for semilinear parabolic equations, cf.~\cite{LUN}, there exist $T_{\max} \in (0,\infty]$ and a unique, maximally defined, classical solution 
$$u \in C\big([0,T_{\max}),X\big) \cap  C^1\big((0,T_{\max}),X\big),$$ of~\eqref{e:intro:main-eq} with initial condition $u(0) = u_0$. If $T_{\max} < \infty$, then it holds $\lim_{t \uparrow T_{\max}} \|u(t)\|_X = \infty$. Proposition~\ref{p:normal-form} then yields a classical solution
$$w \in C\big([0,T_{\max}),X\big) \cap  C^1\big((0,T_{\max}),X\big),$$
of~\eqref{e:weq} given by~\eqref{e:normal-form:def-w}. If $T_{\max} < \infty$, then it holds $\lim_{t \uparrow T_{\max}} \|w(t)\|_X = \infty$. 

The template function $\Theta_u \colon [0,T_{\max}) \to \R$ given by
\begin{equation*}
\Theta_u(t) = \sup_{\mathllap{s} \in [0,t]} \left(\norm{\hat{u}(s)}_{L^\infty} + (1+s)^\frac{1}{2} \norm{\hat{u}(s)}_{L^1}\right),
\end{equation*}
is continuous, monotonically increasing and satisfies $\lim_{t \uparrow T_{\max}} \Theta_u(t) = \infty$ if $T_{\max} < \infty$. Analogously, we define the template function $\Theta_w(t)$ for $w(t)$.

With the notations of Proposition~\ref{p:normal-form}, the normal form transform ensures that 
\begin{equation}
\label{e:nonlinear:normal-form}
\partial_t \hat{w} - \hat{\lambda} \hat{w} = \Cal{F}(Q(u)).
\end{equation}
Using Proposition~\ref{p:normal-form} we bound the right-hand side of~\eqref{e:normal-form:def-w} as
\begin{align*}
\label{e:non-linear:change-of-variable}
\norm{\hat{u}}_{L^p} & {} \leq \norm{\hat{w}}_{L^p} + \norm{\Cal{F}(A_2(u,u))}_{L^p} + \norm{\Cal{F}(A_3(u,u,u))}_{L^p},\\
& {} \lesssim \norm{\hat{w}}_{L^p} + \norm{\hat{u}}_{L^p} \left(\norm{\hat{u}}_{L^1} + \norm{\hat{u}}_{L^1}^2 \right).
\end{align*}
Hence, we obtain
\begin{equation}
\label{e:bound-thetau-thetaw}
\Theta_u(t) \lesssim \Theta_w(t) + \Theta_u(t)^2,
\end{equation}
for $t\in [0,T_{\max})$ with $\Theta_u(t) \leq 1$. Using Proposition~\ref{p:normal-form} again, we estimate the quartic term in~\eqref{e:nonlinear:normal-form} as
\begin{align*}
\norm{\Cal{F}(Q(u))}_{L^p \cap L^\infty} & \lesssim \norm{\Cal{F}(Q(u))}_{L^p} + \norm{\Cal{F}(Q(u))}_{L^\infty}, \\
& {} \lesssim \norm{\hat{u}}_{L^1}^3 \left(\norm{\hat{u}}_{L^p}
+ \norm{\hat{u}}_{L^\infty}\right).
\end{align*}
So, we establish
\begin{equation}
\label{e:nonlinear:bound-Q}
\norm{\Cal{F}(Q(u(s)))}_{L^p \cap L^\infty} \lesssim \frac{\Theta_u(s)^2}{(1+s)^{\frac{3}{2}}},
\end{equation}
for $s \in [0,t]$ with $t\in [0,T_{\max})$ satisfying $\Theta_u(t) \leq 1$. 

We now bound the Duhamel formula associated with~\eqref{e:nonlinear:normal-form}, which reads
\begin{equation*}
\hat{w}(t,k) = \re^{t\hat{\lambda}(k)} \hat{w}(0,k) + \int_0^t \re^{(t-s)\hat{\lambda}(k)} \Cal{F}(Q(u(s)))(k) \dd s.
\end{equation*}
Using the linear bound~\eqref{e:linear-bound} from Lemma~\ref{l:linear-bounds} and the nonlinear bound~\eqref{e:nonlinear:bound-Q}, we get
\begin{align} \label{e:nonlinear-scheme}
\begin{split}
\norm{\hat{w}(t)}_{L^p} \lesssim {} & \frac{\norm{\hat{w}(0)}_{L^p\cap L^\infty}}{(1+t)^\frac{1}{2p}} + \int_0^t \frac{\norm{\Cal{F}(Q(u(s)))}_{L^p\cap L^\infty}}{(1+t-s)^\frac{1}{2p}}\dd s, \\
\lesssim {} & \frac{\norm{\hat{w}(0)}_{L^p\cap L^\infty}}{(1+t)^\frac{1}{2p}} + \Theta_u(t)^2 \int_0^t \frac{\dd s}{(1+t-s)^\frac{1}{2p}(1+s)^{\frac{3}{2}}}, \\
\lesssim {} &  \frac{\norm{\hat{w}(0)}_{L^p\cap L^\infty}}{(1+t)^\frac{1}{2p}} + \frac{\Theta_u(t)^2}{(1+t)^\frac{1}{2p}},
\end{split}
\end{align}
for all $t\in [0,T_{\max})$ with $\Theta_u(t) \leq 1$.

Considering estimate~\eqref{e:nonlinear-scheme} for both $p=1$ and $p=\infty$, and taking the supremum, we conclude that
\begin{equation*}
\Theta_w(t) \lesssim \norm{\hat{w}(0)}_{L^1\cap L^\infty} + \Theta_u(t)^2.
\end{equation*}
holds for all $t\in [0,T_{\max})$ with $\Theta_u(t) \leq 1$. Combining the latter with estimate~\eqref{e:bound-thetau-thetaw}, and applying the estimates~\eqref{e:A2est} and~\eqref{e:A3est} in Proposition~\ref{p:normal-form} at time $t = 0$, we get
\begin{equation*}
\Theta_u(t) \lesssim \norm{\hat{w}(0)}_{L^1\cap L^\infty} + \Theta_u(t)^2 \lesssim \norm{\hat{u}_0}_{L^1\cap L^\infty} + \Theta_u(t)^2,
\end{equation*}
implying that there exists a constant $C > 0$ such that
\begin{equation}
\Theta_u(t) \leq C \left(E_0 + \Theta_u(t)^2\right), 
\end{equation}
for $t\in [0,T_{\max})$ with $\Theta_u(t) \leq 1$, so that taking $E_0 \leq \frac{1}{4C^2}$ yields $T_{\max} = \infty$ and $\Theta_u(t) \leq 2CE_0$ for all $t \geq 0$ by continuity of $\Theta_u$. By interpolation, we deduce that
\begin{equation*}
\norm{\hat{u}(t)}_{L^{1}} \lesssim \frac{\norm{\hat{u}_0}_{L^1\cap L^\infty}}{(1+t)^{\frac{1}{2}}}, \hspace{4em}
\norm{\hat{u}(t)}_{L^{2}} \lesssim \frac{\norm{\hat{u}_0}_{L^1\cap L^\infty}}{(1+t)^{\frac{1}{4}}},
\end{equation*}
\begin{equation*}
\norm{\hat{u}(t)}_{L^{\infty}} \lesssim \norm{\hat{u}_0}_{L^1\cap L^\infty}.
\end{equation*}
holds for all $t \geq 0$. We can then apply the Hausdorff-Young inequality for $p=2,\infty$, and obtain the claimed estimates, which concludes the proof.
\end{proof}

\section{Relation to the space-time resonances method: integration by parts in time} \label{sec:spacetime}

\subsection{Set-up}

The following arguments strongly rely on the properties of the phase function $\phi \colon \R^2 \to \C$ given by 
\begin{equation*}
\phi(k,l) \deq -\hat{\lambda}(k) + \hat{\lambda}(k-l) + \hat{\lambda}(l).
\end{equation*}
We decompose it as $\phi = \rphi + \ri \iphi$, with 
\begin{equation*}
\rphi(k,l) \deq \diffusion \left(k^2 - (k-l)^2 - l^2\right),\\	
\hspace{4em}
\iphi(k,l) \deq \scalp{k} - \scalp{k-l} - \scalp{l}.
\end{equation*}

For the sake of simplicity of exposition, we restrict here to the case $B(u,u) = u^2$. The arguments for general $\hat{B} \in L^\infty(\R^2)$ in~\eqref{e:defB} follow analogously. We apply the Fourier transform to~\eqref{e:intro:main-eq}, and arrive at
\begin{align} \label{e:vHKG2}
\partial_t \hat{u}(t,k) = \hat{\lambda}(k) \hat{u}(t,k) + \hat{u}^{*2}(t,k).
\end{align}
Upon introducing the new coordinate 
\begin{align}\label{e:defv}
v(t,k) = \re^{\ri \scalp{k} t} \hat{u}(t,k),
\end{align}
equation~\eqref{e:vHKG2} transforms into
\begin{align} \label{e:vHKG3}
\partial_t v(t,k) = -\diffusion k^2 v(t,k) + \int_{\R} \re^{\iphi(k,y) \ri t} v(t,k-y)v(t,y) \de y.
\end{align}
The associated Duhamel formulation reads
\begin{align} \label{e:Duhamel}
v(t,k) = \re^{-\diffusion k^2 t} \hat{u}_0(k) + \int_0^t \int_{\R} \re^{-\diffusion k^2(t-s) + \iphi(k,y) \ri s} v(s,k-y)v(s,y) \de y \de s.
\end{align}

Our next step is to transform the quadratic integral term in~\eqref{e:Duhamel} into a quartic one using integration by parts in time. Here, we exploit the fact that the phase function $\phi$ is uniformly bounded away from $0$ on $\R^2$, cf.~Appendix~\ref{app}. In the language of the space-time resonances method of Germain, Masmoudi and Shatah this means that we have an absence of time resonances, cf.~\cite{GERMAIN,LAN,SHAT2}.

\subsection{Integration by parts in time}
Since $\iphi(0,y) \leq -1$ for all $y \in \R$, the mapping $\psi_1 \colon \R^2 \to \R$ given by $\psi_1(k,y) = \diffusion k^2 + \iphi(k,y)\ri$ never vanishes. Thus, for any $f \colon \R^2 \to \C$ satisfying $f(k,y) = f(k,k-y)$ for $(k,y) \in \R^2$, we can integrate by parts in time to establish
\begin{align} \label{e:ibp}
\begin{split}
&\int_0^t \int_{\R} \re^{-\diffusion k^2(t-s) + \iphi(k,y) \ri s} v(s,k-y)v(s,y) f(k,y) \de y \de s\\ 
&\qquad = \left[\int_{\R} \frac{\re^{-\diffusion k^2(t-s) + \iphi(k,y) \ri s}}{\diffusion k^2 + \iphi(k,y)\ri} v(s,k-y)v(s,y) f(k,y) \de y\right]_{s = 0}^t\\[1ex]
&\qquad \quad \ {} -  2\int_0^t \int_{\R^2} \frac{\re^{-\diffusion k^2(t-s) + \left(\iphi(k,y) + \iphi(y,z) \right)\ri s}}{\diffusion k^2 + \iphi(k,y)\ri} v(s,k-y)\\
& \qquad \qquad \qquad \qquad \qquad \qquad \phantom{\int_0^t} \cdot v(s,y-z)v(s,z) f(k,y) \de z\de y \de s \\
&\qquad \quad \ {} + \int_0^t \int_{\R} \frac{\re^{-\diffusion k^2(t-s) + \iphi(k,y) \ri s}}{\diffusion k^2 + \iphi(k,y)\ri} \diffusion \left(y^2 +  (k-y)^2\right)v(s,k-y)v(s,y) f(k,y) \de y \de s,
\end{split}
\end{align}
where we used the equation~\eqref{e:vHKG3} to express the temporal derivatives $\partial_s v(s,y)$ and $\partial_s v(s,k-y)$, and the identity $\iphi(k,y) = \iphi(k,k-y)$ for $(k,y) \in \R^2$ to perform the substition of variables $y \mapsto k-y$ in one of the integrals. 

From Proposition~\ref{p:phi-study:main-result}, the function $\phi$ is bounded away from $0$, uniformly with respect to $(k,y)\in \R^2$. Thus, the function $f_1 \colon \R^2 \to \C$ given by
$$f_1(k,y) = \left(1 - \frac{\diffusion y^2 + \diffusion (k-y)^2}{\diffusion k^2 + \iphi(k,y) \ri}\right)^{-1} = \frac{\diffusion k^2 + \iphi(k,y) \ri}{\phi(k,y)},$$
is well-defined and satisfies $f_1(k,y) = f_1(k,k-y)$ for all $(k,y) \in \R^2$. Thus, subsituting this choice of $f$ into~\eqref{e:ibp} and moving the last integral on the right-hand side of~\eqref{e:ibp} to the left-hand side, we arrive at the identity
\begin{align} \label{e:ibp2}
\begin{split}
&\int_0^t \int_{\R} \re^{-\diffusion k^2(t-s) + \iphi(k,y) \ri s} v(s,k-y)v(s,y) \de y \de s\\ 
&\quad = \left[\int_{\R} \frac{\re^{-\diffusion k^2(t-s) + \iphi(k,y) \ri s}}{\phi(k,y)} v(s,k-y)v(s,y) \de y\right]_{s = 0}^t\\
&\quad \qquad \ {} - 2\int_0^t \int_{\R^2} \frac{\re^{-\diffusion k^2(t-s) + \left(\iphi(k,y) + \iphi(y,z) \right)\ri s}}{\phi(k,y)} v(s,k-y)v(s,y-z)v(s,z) \de z\de y \de s.
\end{split}
\end{align}

Next, we repeat the procedure of integration by parts in time to transform the cubic integral term in~\eqref{e:ibp2} into a quartic one. Since $\iphi(0,y) + \iphi(y,z) \leq -2$ for all $(y,z) \in \R^2$, the mapping $\psi_2 \colon \R^3 \to \R$ given by $\psi_2(k,y,z) = \diffusion k^2 + \iphi(k,y)\ri + \iphi(y,z) \ri$ never vanishes. Thus, for any $g \colon \R^3 \to \C$ and $f \colon \R^2 \to \C$ satisfying $f(k,y) = f(k,k-y)$, $g(k,y,z) = g(k,y,y-z)$ and $g(k,y,z) = g(k,k-z,y-z)$ for $(k,y,z) \in \R^3$, we can integrate by parts in time to obtain
\begin{align} \label{e:ibp3}
\begin{split}
&\int_0^t \int_{\R^2} \re^{-\diffusion k^2t +\psi_2(k,y,z)s} v(s,k-y)v(s,y-z)v(s,z) f(k,y) g(k,y,z) \de z\de y \de s\\ 
&\quad = \left[\int_{\R^2} \frac{\re^{-\diffusion k^2t +\psi_2(k,y,z)s}}{\psi_2(k,y,z)} v(s,k-y)v(s,y-z)v(s,z) f(k,y) g(k,y,z) \de z\de y\right]_{s = 0}^t\\
&\qquad \ {} - \int_0^t \int_{\R^3} \frac{\re^{-\diffusion k^2(t-s) + \left(\iphi(k,y) + \iphi(y,z) + \iphi(z,w) \right)\ri s}}{\psi_2(k,y,z)} v(s,k-y)v(s,y-z)\\
&\qquad \qquad \quad \ \phantom{\int_0^t} \cdot v(s,z-w) v(s,w)\left(2f(k,y) + f(k,y-z)\right) g(k,y,z) \de w\de z\de y \de s\\
&\qquad \ {} + \int_0^t \int_{\R^2} \frac{\re^{-\diffusion k^2t +\psi_2(k,y,z)s}}{\psi_2(k,y,z)} \diffusion \left(z^2 + (y-z)^2 + (k-y)^2\right)\\
&\qquad \qquad \quad \ \phantom{\int_0^t} \cdot v(s,k-y)v(s,y-z)v(s,z) f(k,y) g(k,y,z) \de z\de y \de s,
\end{split}
\end{align}
where we used the equation~\eqref{e:vHKG3} to express the temporal derivatives $\partial_s v(s,k-y)$ and $\partial_s v(s,y-z)$ and $\partial_s v(s,z)$, the identity $\iphi(y,z) = \iphi(y,y-z)$ for $(y,z) \in \R^2$ to substitute the variable $z$ by $y-z$ in one of the integrals, and the identity $\iphi(k,z) + \iphi(k-z,y-z) = \iphi(k,y) + \iphi(y,z)$ for any $(k,y,z) \in \R^3$ for the substitution of variables $(y,z) \mapsto (k-z,y-z)$ in one of the other integrals. 

Proposition~\ref{p:phi-study:main-result} yields that the function $\psi_3 \colon \R^3 \to \C$ given by $\psi_3(k,y,z) = \phi(k,y) + \phi(y,z)$ is bounded away from zero, uniformly with respect to $(k,y,z) \in \R^3$. Thus, the function $g_1 \colon \R^3 \to \C$ given by
$$g_1(k,y,z) = \left(1 - \diffusion \frac{z^2 + (y-z)^2 + (k-y)^2}{\psi_2(k,y,z)}\right)^{-1} = \frac{\psi_2(k,y,z)}{\psi_3(k,y,z)},$$
is well-defined and satisfies $g_1(k,y,z) = g_1(k,y,y-z)$ and $g_1(k,y,z) = g_1(k,k-z,y-z)$ for $(k,y,z) \in \R^3$. Thus, subsituting this choice of $g$ into~\eqref{e:ibp3}, setting $f(k,y) = \smash{\frac{-2}{\phi(k,y)}}$ in~\eqref{e:ibp3} and moving the last integral on the right-hand side of~\eqref{e:ibp3} to the left-hand side we arrive at the identity
\begin{align} \label{e:ibp4}
\begin{split}
&-2\int_0^t \int_{\R^2} \frac{\re^{-\diffusion k^2(t-s) + \left(\iphi(k,y) + \iphi(y,z)\right) \ri s}}{\phi(k,y)} v(s,k-y)v(s,y-z)v(s,z) \de z\de y \de s\\ 
&\qquad = -2\left[\int_{\R^2} \frac{\re^{-\diffusion k^2(t-s) + \left(\iphi(k,y) + \iphi(y,z)\right) \ri s}}{\phi(k,y)\psi_3(k,y,z)} v(s,k-y)v(s,y-z)v(s,z) \de z\de y\right]_{s = 0}^t\\
&\qquad \phantom{{}={}} + 2 \, \int_0^t \int_{\R^3} \frac{\re^{-\diffusion k^2(t-s) + \left(\iphi(k,y) + \iphi(y,z) + \iphi(z,w) \right)\ri s}}{\psi_3(k,y,z)} \left(\frac{2}{\phi(k,y)} + \frac{1}{\phi(k,y-z)}\right)\\
&\qquad \qquad \qquad \qquad \qquad \qquad \phantom{\int_0^t} \cdot  v(s,k-y)v(s,y-z)v(s,z-w)v(s,w) \de w\de z\de y \de s.
\end{split}
\end{align}

All in all, combining~\eqref{e:ibp2} and~\eqref{e:ibp4}, we can rewrite~\eqref{e:Duhamel} as
\begin{align} \label{e:Duhamel2}
\begin{split}
v(t,k) & {} = \re^{-\diffusion k^2 t} \hat{u}_0(k) - \int_{\R} \frac{\re^{-\diffusion k^2t}}{\phi(k,y)} \hat{u}_0(k-y)\hat{u}_0(y) \de y\\[1ex]
& \phantom{={}} + 2 \int_{\R^2} \frac{\re^{-\diffusion k^2t}  \hat{u}_0(k-y)\hat{u}_0(y-z)\hat{u}_0(z)}{\phi(k,y)\psi_3(k,y,z)} \de z\de y \\[1ex]
& \phantom{={}} + \int_{\R} \frac{\re^{\iphi(k,y) \ri t}}{\phi(k,y)} v(t,k-y)v(t,y) \de y\\[1ex]
& \phantom{={}} - 2\int_{\R^2} \frac{\re^{\left(\iphi(k,y) + \iphi(y,z)\right) \ri t}}{\phi(k,y)\psi_3(k,y,z)} v(t,k-y)v(t,y-z)v(t,z) \de z\de y \\[1ex]
& \phantom{={}} + 2 \int_0^t \int_{\R^3} \frac{\re^{-\diffusion k^2(t-s) + \left(\iphi(k,y) + \iphi(y,z) + \iphi(z,w) \right)\ri s}}{\psi_3(k,y,z)} \left(\frac{2}{\phi(k,y)} + \frac{1}{\phi(k,y-z)}\right)\\
&\qquad \qquad \qquad \qquad \qquad \phantom{\int_0^t} \cdot  v(s,k-y)v(s,y-z)v(s,z-w)v(s,w) \de w\de z\de y \de s.
\end{split}
\end{align}

\subsection{Closing the nonlinear iteration}

We close a nonlinear argument by applying standard iterative $L^1$-$L^\infty$-estimates to the Duhamel formulation~\eqref{e:Duhamel2} associated with $v$, which immediately yields the proof of Theorem~\ref{t:main-result} after undoing the change of variables~\eqref{e:defv}.

\begin{proof}[Proof -- Theorem~\ref{t:main-result}]
Let $u_0 \in X$. By standard local existence theory for semilinear parabolic equations there exist $T_{\max} \in (0,\infty]$ and a unique, maximally defined, classical solution 
$$v \in C\big([0,T_{\max}),L^1(\R) \cap L^\infty(\R)\big) \cap  C^1\big((0,T_{\max}),L^1(\R) \cap L^\infty(\R)\big),$$ of~\eqref{e:Duhamel} or, equivalently, of~\eqref{e:Duhamel2} with initial condition $v(0) = \hat{u}_0$. If $T_{\max} < \infty$, then it holds $\lim_{t \uparrow T_{\max}} \|v(t)\|_{L^1 \cap L^\infty} = \infty$. 

Consequently, the template function $\eta \colon [0,T_{\max}) \to \R$ given by
\begin{align*}
\eta(t) = \sup_{s \in [0,t]} \left(\|v(s)\|_\infty + \|v(s)\|_1 \sqrt{1+s}\right), 
\end{align*}
is well-defined, continuous, monotonically increasing and, if $T_{\max} < \infty$, then it satisfies $\lim_{t \uparrow T_{\max}} \eta(t) = \infty$.

Assume $E_0 = \|\hat{u}_0\|_{L^1 \cap L^\infty} \leq 1$ and let $t \in [0,T_{\max})$ be such that $\eta(t) \leq 1$. Using Young's convolution inequality, the integral identity~\eqref{e:intid} and the fact that $\phi$ and $\psi_3$ are uniformly bounded away from $0$ by Proposition~\ref{p:phi-study:main-result}, we estimate the right-hand side of~\eqref{e:Duhamel2} by
\begin{align*}
\|v(t)\|_\infty \lesssim E_0 + \frac{\eta(t)^2}{\sqrt{1+t}} + \int_0^t \frac{\eta(s)^2}{(1+s)^{\frac{3}{2}}} \de s  \lesssim E_0 + \eta(t)^2,
\end{align*}
and
\begin{align*}
\|v(t)\|_1 \lesssim \frac{E_0}{\sqrt{1+t}} + \frac{\eta(t)^2}{1+t} + \int_0^t \frac{\eta(s)^2}{\sqrt{t-s}(1+s)^{\frac{3}{2}}} \de s  \lesssim \frac{E_0 + \eta(t)^2}{\sqrt{1+t}}.
\end{align*}
Combining the latter two estimates yields a constant $C \geq 1$ such that for all $t \in [0,T_{\max})$ with $\eta(t) \leq 1$ it holds
\begin{align*}
\eta(t) \leq C\left(E_0 + \eta(t)^2\right),
\end{align*}
so that $E_0 \leq \frac{1}{4C^2}$ yields $T_{\max} = \infty$ and $\eta(t) \leq 2CE_0$ for all $t \geq 0$. 

Thus, undoing the coordinate transform implies that there exist constants $K, \delta > 0$ such that for each $u_0 \in X$ satisfying $E_0 := \|u_0\|_X < \delta$ there exists a unique global classical solution 
$$u \in C\big([0,\infty),X\big) \cap C^1\big((0,\infty),X\big),$$ of~\eqref{e:intro:main-eq} with initial condition $u(0) = u_0$ obeying the estimates
\begin{align*}
\|u(t)\|_2 \leq KE_0(1+t)^{-\frac14}, \qquad \|u(t)\|_{\infty} \leq \|\hat{u}(t)\|_1 \leq KE_0(1+t)^{-\frac12}, \qquad \|\hat{u}(t)\|_\infty \leq KE_0,
\end{align*}
for all $t \geq 0$, which completes the proof.
\end{proof}

\appendix

\section{Phase study} \label{app}

We derive lower bounds on the phase function $\phi \colon \R^2 \to \C$ given by 
\begin{equation*}
\phi(k,l) \deq -\hat{\lambda}(k) + \hat{\lambda}(k-l) + \hat{\lambda}(l).
\end{equation*}

\begin{proposition}
\label{p:phi-study:main-result}
There exists a constant $\phi_0 > 0$ such that 
\begin{equation*}
\absolute{\phi(k,l)} \geq \phi_0,
\hspace{4em}
(k,l) \in \R^2,
\end{equation*}
\begin{equation*}
\absolute{\phi(k,l) + \phi(l,m)} \geq \phi_0,
\hspace{4em}
(k,l,m) \in \R^3.
\end{equation*}
\end{proposition}
\begin{proof}
As before, we decompose $\phi$ into its real and imaginary part as $\phi = \rphi + \ri \iphi$ with 
\begin{equation*}
\rphi(k,l) = d\left(k^2 - (k-l)^2 - l^2\right),
\hspace{4em}
\iphi(k,l) =\scalp{k} - \scalp{k-l} - \scalp{l}.
\end{equation*}
First, if $(k,l) \in \R^2$ satisfies $|l| \geq \frac{1}{2}$ and $|k-l| \geq \frac{1}{2}$, then it holds 
$$|\rphi(k,l)| = 2d |l||k-l| \geq \frac{d}{2} > 0.$$ 
On the other hand, if $(k,l) \in \R^2$ satisfies $|l| \leq \frac{1}{2}$, then the convexity of $g \colon \R \to \R, g(k) = \scalp{2k}$ yields a constant $\delta > 0$ such that
\begin{align*} \iphi(k,l) &= g\left(\frac12 l + \frac12 (k-l)\right) - \scalp{l} - \scalp{k-l}\\ &\leq \frac12 g(l) + \frac12 g(k-l) - \scalp{l} - \scalp{k-l} \leq \sqrt{\frac14 + l^2} - \sqrt{1+l^2} < -\delta.\end{align*}
Similarly, we find for $(k,l) \in \R^2$ with $|k-l| \leq \frac{1}{2}$ that $\iphi(k,l) < -\delta$. Moreover, we showed in passing that $\iphi(k,l)$ is negative for any $(k,l) \in \R^2$. We conclude that $\phi$ is uniformly bounded away from $0$ on $\R^2$. 

Subsequently, we consider the mapping $\psi \colon \R^3 \to \C$ given by $\psi(k,l,m) = \phi(k,l) + \phi(l,m)$. If $(k,l,m) \in \R^3$ is such that $\mathrm{sgn}(l) = \mathrm{sgn}(k-l)$, $\mathrm{sgn}(m) = \mathrm{sgn}(l-m)$ and $|l|,|k-l|,|m|,|l-m| \geq \frac{1}{2}$, then it holds 
$$\mathrm{Re}(\psi(k,l,m)) = \rphi(k,l) + \rphi(l,m) = 2d \left(l (k-l) + m(l-m)\right) \geq d > 0.$$ 
On the other hand, if $(k,l) \in \R^2$ is such that $k$ and $k-l$ have opposite sign, then we have $|k| \leq |l|$, implying $\scalp{l} \geq \scalp{k}$. So, we have $\iphi(k,l) \leq -1$. Hence, if $(k,l,m) \in \R^3$ is such that $\mathrm{sgn}(l) \neq \mathrm{sgn}(k-l)$ or $\mathrm{sgn}(m) \neq \mathrm{sgn}(l-m)$, then it holds $$\mathrm{Im}(\psi(k,l,m)) = \iphi(k,l) + \iphi(l,m) \leq -1 < 0,$$ 
using that $\iphi$ is negative on $\R^2$. Finally, if $(k,l,m) \in \R^3$ satisfies $|l| \leq \frac{1}{2}$, $|k-l| \leq \frac{1}{2}$, $|l-m| \leq \frac{1}{2}$ or $|m| \leq \frac{1}{2}$, then we have $$\mathrm{Im}(\psi(k,l,m)) = \iphi(k,l) + \iphi(l,m) \leq -\delta,$$ using that $\iphi$ is negative on $\R^2$, and that $\iphi(k,l) \leq -\delta$ holds for $(k,l) \in \R^2$ with $|k-l| \leq \frac12$ or $|k| \leq \frac12$ (as derived above). We conclude that $\psi$ is uniformly bounded away from $0$ on $\R^3$.
\end{proof}

\bibliographystyle{abbrv}
\bibliography{KGbib}

\end{document}